\documentclass[11pt]{article}

\usepackage{amssymb,amsmath,latexsym}
\usepackage[dvips]{graphics}
\usepackage[mathscr]{eucal}

\def\doublespace{\baselineskip=14pt}

\def\EE{{\cal E}}
\def\FF{{\cal F}}

\def\DD{{\cal D}}
\def\LL{{\cal L}}
\def\NN{{\cal N}}
\def\AA{{\cal A}}

\def\hh{{\cal H}}

\def\E{{\bf E}}

\def\P{{\bf P}}
\def\Q{{\bf Q}}
\def\H{{\bf H}}

\def\RR{\mathbb{R}}
\def\PP{\mathbb{P}}

\def\HH{\mathbb{H}}

\def\a{{\bf a}}

\def\n{{\bf n}}
\def\D{{\bf D}}

\def\1{{\bf 1}}
\def\<{\langle}
\def\>{\rangle}
\def\pf{\noindent{\bf Proof.} }

\def\wt{\widetilde}
\def\wh{\widehat}
\def\dis{\displaystyle}

\def\o{{\bf 0}}

\def\qed{{\hfill $\Box$ \bigskip}}

\numberwithin{equation}{section}

\newtheorem{thm}{Theorem}[section]
\newtheorem{prop}[thm]{Proposition}
\newtheorem{remark}[thm]{Remark}
\newtheorem{defn}[thm]{Definition}

\begin{document}

\title
{\bf Reflections at infinity of time changed RBMs \\
on a domain with Liouville branches} 
\author{\ Zhen-Qing Chen and Masatoshi Fukushima}
\date{October 3, 2016}
\maketitle

\begin{abstract}
Let $Z$ be the transient reflecting Brownian motion on the closure of an unbounded domain $D\subset \RR^d$ with $N$ number of Liouville branches.  
We consider a diffusion 
$X$
on $\overline D$ having finite lifetime obtained from $Z$ by a time change.  We show that $X$ admits only a finite number of possible symmetric conservative diffusion extensions $Y$ beyond its lifetime 
characterized by possible partitions of the collection of $N$ ends
and we identify the family of the extended Dirichlet spaces of all $Y$ (which are independent  of  time change used)
as subspaces of the space ${\rm BL}(D)$ 
spanned by the extended Sobolev space $H_e^1(D)$ and the approaching probabilities of $Z$ to the ends of Liouville branches.
\end{abstract}

\medskip
\noindent
{\bf AMS 2010 mathematics Subject Classification}: Primary 60J50, Secondary 60J65, 32C25 


 \doublespace

\section{Introduction}\label{S:1}
The boundary problem of a Markov process $X$ concerns all possible Markovian prolongations $Y$ of $X$ beyond its life time $\zeta$ whenever $\zeta$ is finite.   For a conservative but transient Markov process, we can still consider its extension,  after a time change to speed up the original process.
Let $Z=(Z_t, \Q_z)$ be a conservative right process on a locally compact separable metric space $E$ and $\partial$ be the point at infinity of $E$.  Suppose $Z$ is transient relative to an excessive measure $m$: for the $0$-order resolvent $R$ of $Z$, $Rf(z)<\infty,\ m$-a.e. for some strictly positive function (or equivalently, for any non-negative function) $f\in L^1(E;m)$.  
Then
\[\Q_z\left(\lim_{t\to\infty}\;Z_t\;=\;\partial\right)=1\quad \text{for q.e.}\ x\in E,\]
if $Rf$ is lower semicontinuous for any non-negative Borel function $f$ (\cite{FTa}).  The last condition is not needed when $X$ is $m$-symmetric (\cite{CF2}).
Here, 'q.e.' means 'except for an $m$-polar set.

Take any strictly positive bounded function $f\in L^1(E;m).$   Then $A_t=\int_0^t f(Z_s)ds,\ t\ge 0$ is a strictly increasing PCAF of $Z$ with $\E_z^\Q[A_\infty]=Rf(x)<\infty$ for q.e. $x\in E.$ \\
The time changed process $X=(X_t, \zeta,\P_x)$ of $Z$ by means of $A$ is defined by
\[X_t=Z_{\tau_t},\ t\ge 0,\quad \tau=A^{-1},\quad \zeta=A_\infty,\quad \P_x=\Q_x,\ x\in E.\]
Since
$\dis \P_x(\zeta<\infty,\ \lim_{t\to\zeta}X_t=\partial)=\P_x(\zeta<\infty)=1$ for q.e. $x\in E,$
the boundary problem for $X$ at $\partial$ makes perfect sense.  We denote $X$ also by $X^f$ to indicate its depedence on the function $f$.  For different choices of $f$, $X^f$ have a common geometric structure related each other only by time changes.  Thus a study of the boundary problem for $X=X^f$ is a good way to have a closer look 
at the 
geometric behaviors of a conservative transient process $Z$ around $\partial.$ 
A strong Markov process $\wh X$ on a topological space $\wh E$ is said to be an extension of $X$ on $E$ if
(i) $E$ can be embedded  homomorphically as a dense open subset of $\wh E$, (ii) the part process of $\wh X$ killed 
upon leaving $E$ has the same distribution as $X$, and (iii) $\wh X$ has no sojourn on $\wh E\setminus E$; 
that is, $\wh X$ spends zero Lebesgue amount of time on $\wh E \setminus E$.  

In this paper, we consider as $Z$ the transient reflecting Brownian motion on the closure of an unbounded domain $D\subset \RR^d$ with $N$ number of Liouville branches.  Our main aim is 
to prove in Section \ref{S:4} that a time changed process $X^f$ of $Z$ 
admits essentially 
only a finite number of possible symmetric conservative diffusion extensions $Y$ beyond its 
lifetime. They are characterized 
by the partition of the collection of $N$ ends. 
Moreover,  all the corresponding extended Dirichlet spaces $(\EE^Y, \FF^Y_e)$ are identified 
in terms of the extended Dirichlet space of $Z$ and the approaching probabilities of $Z$ to the ends of Liouville branches in an extremely simple manner.
These extended Dirichlet spaces are 
independent of the choice of $f.$  The $L^2$-generator of each extension $Y$ 
is 
also characterized in Section 
\ref{S:5} by means of 
zero flux conditions at  the ends of branches.
Each extension $Y$ may be called a {\it many point reflection at infinity of} $X^f$ generalizing the notion of the one point reflection in \cite{CF3} in the present specific context.
The characterization of possible extensions also uses  quasi-homeomorphism and equivalence between Dirichlet forms. 
See the Appendix, Section \ref{S:7}, of this paper for details. 

In fact, our results are valid  
for a time changed process $X^\mu$ of $Z$ by means of a more general finite smooth measure $\mu$ on $\overline D$ than $f (x)dx.$
This is demonstrated in Section \ref{S:6}.  

\smallskip

Although we formulate our results for the reflecting Brownian motion on an unbounded domain in $\RR^d$ with several Liouville branches, all of them except for Theorem \ref{T:5.1} remain valid without any essential change for the reflecting diffusion process associated with the uniformly elliptic second order self-adjoint partial differential operator with measurable coefficients that was constructed in 
\cite{C} and \cite{FTo}. 
Since we need strong Feller property of the reflecting diffusion process, 
we assume the underlying unbounded domain is Lipschitz  in the sense of \cite{FTo}; 
see Remark \ref{R:4.3}. 
Thus we are  effectively 
 investigating common path behaviours at infinity holding for such a general family of diffusion processes. 
 
\smallskip 

{\bf Acknowledgement}\quad This paper is a direct outgrowth of  our paper \cite{CF1} and Chapter 7 of our book \cite{CF2}. 
In relation to them,  
we had very valuable discussions with Krztsztof Burdzy  on boundaries of transient reflecting Brownian motions.  
We would like to express our sincere thanks to him. 

\section{Preliminaries}

For a domain $D\subset \RR^d$, let us consider the spaces 
\begin{equation}\label{e;1.1}
{\rm BL}(D)=\{u\in L^2_{\rm loc}(D):\; |\nabla u|\in L^2(D)\},\quad
 H^1(D)={\rm BL}(D)\cap L^2(D), 
\end{equation} 
The space ${\rm BL}(D)$ called the 
{\it Beppo Levi space} was introduced by J. Deny and J. L. Lions \cite{DL} 
as the space of Schwartz distributions whose first order derivatives are in $L^2(D)$, which can be identified with the function space described above.  The quotient space $\dot{\rm BL}(D)$ of ${\rm BL}(D)$ by the space of all constant functions on $D$ is a real Hilbert space with inner product
\[\D(u,v)=\int_D\nabla u(x)\cdot\nabla v(x)dx.\]
See \S 1.1 of V.G, Maz'ja \cite{M} for proofs of 
the above
 stated facts, 
where the space BL$(D)$ is denoted by $L_2^1(D)$ and studied in a more general context of the spaces $L_p^\ell(D),\; \ell \ge 1,\; p\ge 1.$ 

Define
\begin{equation}\label{e:1.2}
(\EE, \FF)=(\frac12 \D, H^1(D)),
\end{equation} 
which is a Dirichlet form on $L^2(D)$.  
The collection of those domains $D\subset \RR^d$ for which \eqref{e:1.2} is regular on $L^2(\overline D)$ will be denoted by $\DD.$  
It is known that 
$D\in \DD$ if $D$ is either a domain of continuous boundary or an extendable domain relative to $H^1(D)$ (cf. \cite[p 866]{CF1}).   
For $D\in \DD,$
the diffusion process $Z$ on $\overline D$ associated with \eqref{e:1.2} is by definition the {\it reflecting Brownian motion} (RBM in abbreviation) which is known to be conservative. 
 Furthermore, the space ${\rm BL}(D)$ is nothing but the {\it reflected Dirichlet space} of the form \eqref{e:1.2} (\cite[\S 6.5]{CF2}).  
The Dirichlet form \eqref{e:1.2} is either recurrent or transient and the latter case occurs only when $d\ge 3$ and $D$ is unbounded. For $D_1,\; D_2 \in \DD$ with $D_1\subset D_2$, \eqref{e:1.2} is transient for $D_2$ whenever so it is for the smaller domain $D_1$.  If \eqref{e:1.2} is recurrent, then we have the identity
\[ {\rm BL}(D)=H_e^1(D)\]
where $H_e^1(D)$ denotes the {\it extended Dirichlet space} of 
the form \eqref{e:1.2} or of the RBM $Z$ (\cite{CF2}) 
that may be called the {\it extended Sobolev space of order} $1$.

Suppose $D\in \DD$ and 
 \eqref{e:1.2} is transient.  Then $H_e^1(D)$ is a Hilbert space with inner product $\frac12 \D$ possessing the space $C_c^\infty(\overline D)$ as its core.  $H_e^1(D)$ can be regarded as a proper closed subspace of the quotient space $\dot{\rm BL}(D).$  Define
\begin{equation}\label{e:1.3}
\HH^*(D)=\{u\in {\rm BL}(D): \ \D(u,v)=0\ \hbox{\rm for every}\ v\in H^1_e(D)\}.
\end{equation}
Any function $u\in {\rm BL}(D)$ admits a unique decomposition
\begin{equation}\label{e:1.4a}
u = u_0+ h,\quad u_0\in H_e^1(D),\quad h\in \HH^*(D).
\end{equation}
Any function $h\in \HH^*(D)$ is of finite Dirichlet integral and harmonic on $D.$  Furthermore, the quasi-continuous version of $h$ is harmonic on $\overline D$ with respect to the RBM $Z$.
 
In what follows, we restrict our attention to the case where the form \eqref{e:1.2} is transient and so we assume that $d\ge 3$ and $D\in \DD$ is unbounded.

\begin{defn}\label{D:1.1} \rm A domain $D\in \DD$ is called a {\it Liouville domain} if
the form \eqref{e:1.2} is transient and ${\rm dim} \HH^*(D)=1.$
\end{defn}

A domain $D\in \DD$ is a Liouville domain if and only if the form \eqref{e:1.2} is transient and
any function $u\in{\rm BL}(D)$ admits a unique decomposition
\begin{equation}\label{e:1.4}
u = u_0 +c, \qquad \hbox{where } u_0\in H_e^1(D) \hbox{ and } c\in \RR. 
\end{equation}
We shall denote by $c(u)$ the constant $c$ in \eqref{e:1.4} uniquely associated with $u\in {\rm BL}(D)$ for a Liouville domain $D$.
 
A trivial but important example of a Liouville domain is  $\RR^d$ with $ d\ge 3$,  see M. Brelot \cite{B}.
Another important example of a Liouville domain is provided by an unbounded uniform domain that has been shown by  
P. Jones \cite{J} (see also \cite{HK}) 
to be an extendable domain relative  to the space ${\rm BL}(D)$.

A domain $D\subset \RR^d$ is called a {\it uniform domain} if there exists $C>0$ such that for every $x,y\in D$, 
 there is a rectifiable curve $\gamma$ in $D$ connecting $x$ and $y$ with length$(\gamma)\le C|x-y|$, and moreover
\[\min\{ |x-z|,\;|z-y|\}\le C{\rm dist}(z,D^c)\quad \hbox{\rm for every}\ z\in \gamma.\]
It was proved in Theorem 3.5 of \cite{CF1} that any unbounded uniform domain is a Liouville domain in the sense of Definition \ref{D:1.1}.  An unbounded uniform domain is such a domain that is broaden toward the infinity.  The truncated infnite cone $C_{A,a}=\{(r,\omega): r>a, \;\omega\in A\}\subset \RR^d$ for any connected open set $A\subset S^{d-1}$ with Lipschitz boundary is an unbounded uniform domain.  
To the contrary, \eqref{e:1.2} is recurrent for the cylinder $D=\{(x,x')\in \RR^d: x\in \RR,\ |x'|<1\}$.  
See R. G. Pinsky \cite{P} 
for transience criteria for other types of domains.  
On the other hand, it has been shown in 
\cite[Proposition 7.8.5]{CF2} 
that \eqref{e:1.2} is transient but ${\rm dim}(\HH^*(D))=2$ for a special domain 
\begin{equation}\label{e:1.5}
D=B_1(\o)\:\cup\; \left\{(x,x')\in \RR^d: x\in \RR,\ |x|> |x'|\right\},\quad d\ge 3.
\end{equation}
with two symmetric cone branches.  Here $B_r(\o),\;r>0,$ denotes an open ball with radius $r$ centered at the origin.  This domain is not uniform because of a presence of a bottleneck.  We shall consider much more general domains than this. But before proceeding to the main setting of the present paper,  
we state a simple property of Liouville domains:

\begin{prop}\label{P:1.2}
For $D_1, D_2\in \DD$ with $D_1\subset D_2,$ Suppose $D_1$ is a Liouville domain and $D_2\setminus D_1$ is bounded.   
Then $D_2$ is a Liouville domain.  
Furthermore, for any $u\in {\rm BL}(D_2)$, it holds that $c(u)=c(u\big|_{D_1}).$
\end{prop}

\pf The proof is similar to that of \cite[Proposition 3.6]{CF1}. 
 Note that \eqref{e:1.2} is transient for $D_2.$  We show that any $u\in {\rm BL}(D_2)$ admits a decomposition \eqref{e:1.4} with $u_0\in H_e^1(D_2)$ and $c=c(u\big|_{D_1})$. 
 Due to the normal contraction property of ${\rm BL}(D_2)$ and the transience of $( \frac12 {\bf D}, H^1(D))$, 
 we may assume that $u$ is bounded on $D_2.$  By noting that $u\big|_{D_1}\in {\rm BL}(D_1)$ and $D_1$ is a Liouville domain, we let $c=c(u\big|_{D_1})$ and $u_0(x)=u(x)-c,\ x\in D_2.$  Then $u_0\big|_{D_1}\in H_e^1(D_1).$  To prove that $u_0\in H_e^1(D_2),$ choose an open ball $B_r(\o)\supset \overline{D_2\setminus D_1}$ and a function $w\in C_c^\infty(\RR^d)$ with $w(x)=1,\ x\in B_r(\o).$  Clearly $wu_0\in H_e^1(D_2).$  

It remains to show $(1-w)u_0\in H_e^1(D_2).$  Take $g_n\in H^1(D_1)$ converging to $u_0$ a.e. on $D_1$ and in the Dirichlet norm on $D_1$.  By truncation, we may assume that $g_n$ is uniformly bounded on $D_1.$  Then
\begin{eqnarray*} 
&&\int_{D_2} |\nabla [(1-w(x))g_n(x)]|^2dx\\
& \le & 2 \sup_{x\in \RR^d} (1-w(x))^2 \int_{D_1}|\nabla g_n(x)|^2dx + 2\sup_{x\in D_1} |g_n(x)|^2\int_{\RR^d} |\nabla w(x)|^2 dx,
\end{eqnarray*}
which is uniformly bounded in $n$, yielding by the Banach-Saks theorem that $(1-w)u_0\in H_e^1(D_2).$ \qed

We shall work under the regularity condition

\smallskip\noindent
{\bf(A.1)}\ $D$ is of a Lipschitz boundary $\partial D,$

\smallskip\noindent
which means the following: there are constants $M>0,\; \delta>0$ and a locally finite covering $\{U_j\}_{j\in J}$ of $\partial D$ such that, for each $j\in J,$ $D\cap U_j$ is a upper part of a graph of a Lipschitz continuous function under an appropriate coodinate system with the Lipschitz constant bounded by $M$ and $\partial D\subset \bigcup_{j\in J}\{x\in U_j: {\rm dist}(x,\partial U_j)>\delta\}.$   
According to \cite{FTo}, there exists then a conservative diffusion process $Z=(Z_t, \Q_x)$ on $\overline D$ associated with the regular Dirichlet form \eqref{e:1.2} on $L^2(\overline D)$ whose resolvent $\{G_\alpha^Z; \alpha>0\}$ has the strong Feller property in the sense that
\begin{equation}\label{e:1.6}
G_\alpha^Z(bL^1(D))\subset bC(\overline D).
\end{equation} 
$Z$ is a precise version of the RBM on $\overline D.$  In particular, the transition probability of $Z$ is absolutely continuous with respect to the Lebesgue measure.  

Under the condition {\bf (A.1)} and the transience assumption on \eqref{e:1.2}, the RBM $Z=(Z_t, \Q_x)$ on $\overline D$ enjoys the properties that
\begin{equation}\label{e:1.7} 
\Q_x\left(\lim_{t\to\infty} Z_t=\partial\right)=1
\quad\text{for every}\ x\in \overline D,
\end{equation}
where $\partial$ denotes the point at infinity of $\RR^d$, and
\begin{equation}\label{e:1.8}
\Q_x\left(\lim_{t\to\infty} u(Z_t)=0\right)=1
\quad\text{for every}\ x\in \overline D,
\end{equation}
for any $u\in H_e^1(D)$, $u$ being taken to be quasi-continuous.  
See \cite[\S 7.8, $(4^o)$]{CF2}. 
 
In the rest of this paper, we fix a domain $D$ of $\RR^d,\; d\ge 3,$ satisfying {\bf(A.1)} and

\smallskip\noindent  
{\bf (A.2)}\qquad $\dis D\;\setminus\;\overline{B_r(\o)}\;=\;\bigcup_{j=1}^N\; C_j$

\smallskip\noindent  
for some $r>0$ and an integer $N$,\ where $C_1,\cdots, C_N$ are  
Liouville domains with Lipschitz boundaries such that $\overline C_1, \cdots, \overline C_N$ are mutually disjoint.  
$D$ may be called a {\it Lipschitz domain with $N$ number of Liouville branches}.    

Let $\partial_j$ be the point at infinity of the unbounded closed set $\overline C_j$ for each $1\le j\le N.$
Denote the $N$-points set $\{\partial_1,\cdots, \partial_N\}$ by $F$ and put $\overline D^*=\overline D\cup F.$  $\overline D^*$ can be made to be a compact Hausdorff space if we employ as a local base of neighborhoods of each point $\partial_j\in F$ the neighborhoods of $\partial_j$ in $\overline C_j\cup \{\partial_j\}.$  $\overline D^*$ may be called the {\it $N$-points compactification of} $\overline D.$

Obviously the Dirichlet form \eqref{e:1.2} is transient for $D$.   We shall verify in Section \ref{S:3} that dim$(\HH^*(D))=N.$   
Here we note the following implication of Proposition \ref{P:1.2}; if a domain $D$ is of the type {\bf (A.2)} for different $0<r_1<r_2,$ and if $D$ is a domain with $N$ number of Liouville branches relative to $r_2$, then so it is relative to $r_1.$

\section{Approaching probabilities of RBM $Z$ and limits of BL-functions along $Z_t$}\label{S:2}

For each $1\le j\le N,$ define the approaching probability of the RBM $Z=(Z_t, \Q_x)$ to $\partial_j$ by
\begin{equation}\label{2.1}
\varphi_j(x)=\Q_x\left(\lim_{t\to\infty}Z_t=\partial_j\right),\quad x\in \overline D.
\end{equation}  

\begin{prop}\label{P:2.1}\ It holds that
\begin{equation}\label{2.2}
\sum_{j=1}^N \varphi_j(x)=1\quad\text{for every}\quad x\in \overline D,
\end{equation}
and, for each $1\le j\le N,$
\begin{equation}\label{2.3}
 \varphi_j(x)>0\quad\text{for every}\quad x\in \overline D.
\end{equation}

\end{prop}
\pf \eqref{2.2} is a consequence of \eqref{e:1.7}.  
As $\varphi_j$ is a non-negative harmonic function on the domain $D$, it is either identically zero on $D$ or  
strictly positive on $D$.
Since $\varphi_j(x) = Q_t \varphi_j(x),\ x\in \overline D$, where $Q_t$ is the transition semigroup of the RBM $Z$,
 which has a strictly positive transition density kernel,
 the above dichotomy  extends from $D$ to $\overline D.$ 

Suppose $\varphi_j(x) \equiv 0$ on $\overline D$. 
Then by \eqref{e:1.7}
\begin{equation}\label{2.4}
\Q_x\left(\sigma_{\partial B_r(\o)}<\infty\right)=1,\quad \text{\rm for any}\ x\in \overline C_j\setminus B_{r+1}(\o).
\end{equation}
Let $Z^j=(Z_t^j,\Q_x^j),\ x\in \overline C_j,$ be the RBM on $\overline C_j$,
which is transient as $C_j$ is a Liouville domain. 
Since $Z$ and $Z^j$ share the common part process on $\overline C_j\setminus \partial B_r(\o),$ \eqref{2.4} remains valid if $\Q_x$ is replaced by $\Q_x^j.$  By the Markov property of $Z^j$ and the consevativeness of $Z^j,$
 we have  
\[\Q_x^j\left(\sigma_{\partial B_r(\o)}\circ \theta_\ell<\infty\ \text{\rm for every integer}\ \ell\right)=1,\]
for any $x\in \overline C_j\setminus B_{r+1}(\o)$.
This however contradicts to the transience property \eqref{e:1.7} of $Z^j$. \qed

\begin{prop}\label{P:2.2}
 For any $u\in {\rm BL}(D)$,  
let $
c_j (u)
= c(u|_{C_j})$ for $1\leq j\leq N$. Then
\begin{equation}\label{3.6} 
\Q_x\left(Z_{\infty-}=\partial_j,\ \lim_{t\to\infty} u(Z_t)=c_j (u)\right)
=\Q_x\left(Z_{\infty-}=\partial_j\right),\ x\in \overline D,\quad 1\le j\le N.
\end{equation} 
If $c_j (u)=0$ for every $1\le j\le N,$ then $u\in H_e^1(D).$  
\end{prop}

\pf\  We prove \eqref{3.6} for $j=1.$  
Let $r>0$ be the radius in {\bf(A.2)} and $Z^1=(Z_t^1, \Q_x^1)$ be the RBM
on $\overline {C_1}$.  The hitting times of $B_r(\o)$ and $B_R(\o)$ for $R>r$ will be denoted by $\sigma_r$ and $\sigma_R$, respectively.
Observe that $Z$ and $Z^1$ share in common the part process on $\overline C_1\setminus \partial B_r(\o).$
Since $C_1$ is a Liouville domain, we see from \eqref{e:1.4} and \eqref{e:1.8} that
$$ 
\Q_x^1 \left( \lim_{t\to \infty} u(Z_t^1) = c_1 (u)\right)=1
\quad \hbox{for every } x\in \overline C_1.
$$

For $R>r$, we consider the event
\[ \Gamma_R=\{Z_{\sigma_R}\in \overline C_1,\ \sigma_r\circ \theta_{\sigma_R}=\infty\}.\]
Then $\Gamma_R\cap\{Z_{\infty-}=\partial\}$ increses as $R$ increases and
$\{Z_{\infty-}=\partial_1\}=\bigcup_{R>r} [\Gamma_R\cap \{Z_{\infty-}=\partial\}]$.   
In view of \eqref{e:1.7}, we  have for $x\in \overline D$, 
\begin{eqnarray*}
&& \Q_x(Z_{\infty-}=\partial_1)
= \lim_{R\to\infty} \Q_x(\Gamma_R\cap \{Z_{\infty-}=\partial\})
= \lim_{R\to\infty} \Q_x(\Gamma_R)\\
&&=\lim_{R\to\infty} \E^{\Q_x}\left[\Q_{Z_{\sigma_R}}(\sigma_r=\infty); Z_{\sigma_R}\in \overline C_1\right]\\
&&=\lim_{R\to\infty} \E^{\Q_x}\left[\Q_{Z_{\sigma_R}}^1(\sigma_r=\infty); Z_{\sigma_R}\in \overline C_1\right]\\
&&=\lim_{R\to\infty} \E^{\Q_x}\left[\Q_{Z_{\sigma_R}}^1(\sigma_r=\infty,\;\lim_{t\to\infty}  
u(Z_t^1)=c_1 (u)); Z_{\sigma_R}\in \overline C_1\right] . 
\end{eqnarray*}
In exactly the same way, we can see that  
$\Q_x(Z_{\infty-}=\partial_1,\; \lim_{t\to\infty} u(Z_t)=c_1 (u))$ 
equals the last expression in the above display, 
 proving \eqref{3.6} for $j=1.$

Suppose $u\in {\rm BL}(D)$ satisfies $c_j (u)=0$ for every $1\le j\le N.$ 
Then $u\big|_{C_j} \in H^1_e(C_j)$ for every $1\le j\le N$ and we can conclude as the proof of Proposition
\ref{P:1.2} that $u\in H_e^1(D).$  \qed

We remark that, in view of Proposition \ref{P:1.2} the constants  $c_j (u),\:1\le j\le N,$ 
in the above proposition are independent of the choice of the radius $r$ in {\bf (A.2)}.

\section{Reflecting extension $X^*$ of a time changed RBM $X$ and dimension of $\HH^*(D)$}
\label{S:3}

Fix a strictly positive bounded integrable function $f$ on $\overline D$ and define
\begin{equation}\label{e:3.1}
A_t=\int_0^tf(Z_s)ds,\quad t\ge 0.
\end{equation}
$A_t$ is a positive continuous additive functional (PCAF) of the RBM $Z=(Z_t, \Q_x)$ on $\overline D$ in the strict sense with full support.  Notice that
\begin{equation}\label{3.1}
\Q_x(A_\infty<\infty)=1\quad \text{for every}\ x\in \overline D,
\end{equation} 
because $\E^{Q_x}\left[A_\infty\right]=G_{0+}^Zf(x)<\infty$ for a.e. $x\in \overline D$ due to the transience of $Z$ (\cite[Proposition 2.1.3]{CF2}) and hence
\begin{equation} \label{e:3.2} 
\Q_x(A_\infty=\infty)=\Q_x(A_\infty\circ \theta_t=\infty)=\E^{\Q_x}\left[\Q_{Z_t}(A_\infty=\infty)\right]=0
\quad \hbox{for every }  x\in \overline D,
\end{equation}
on account of the stated absolute continuity of the transition function of $Z.$

Let $X=(X_t, \zeta, \P_x)$ be the time changed process of $Z$ by means of $A$:
\[X_t=Z_{\tau_t},\quad \tau=A^{-1},\quad \zeta=A_\infty,\quad \P_x=\Q_x\ {\rm for}\ x\in \overline D.\] 
The Markov process 
$X=X^f$ 
is a diffusion process on $\overline D$ symmetric with respect to the measure $m(dx)=f(x)dx$ and the Dirichlet form $(\EE^X, \FF^X)$ of $X$ on $L^2(\overline D;m)$ is given by
\begin{equation}\label{3.2}
\EE^X=\frac12\D, \qquad \FF^X=H_e^1(D)\cap L^2(\overline D;m).
\end{equation}
Since the extended Dirichlet space and the reflected Dirichlet space are invariant under a time change by a fully supported PCAF (\cite[Cor.5.2.12,\; Prop.6.4.6]{CF2}), these spaces for $\EE^X$ are still given by $H_e^1(D)$ and BL$(D)$, respectively.  But the life time $\zeta$ of $X$ is finite $\P_x$-a.s. for every $x\in \overline D$ in view of \eqref{3.1} so that we may consider the problem of extending $X$ after $\zeta,$ particularly, from $\overline D$ to its $N$-points compactification $\overline D^*=\overline D\cup F$ with $F=\{\partial_1, \cdots, \partial_N\}.$ 

We can rewrite the approaching probability $\varphi_j$ of $Z$ to $\partial_j$ defined by \eqref{2.1} as
\begin{equation}\label{3.3}
\varphi_j(x)=\P_x\left(\zeta<\infty,\quad X_{\zeta-}=\partial_j\right),\quad x\in \overline D,\quad 1\le j\le N,
\end{equation}
in terms of the time changed process $X$.   The measure 
$m(dx)=f(x)dx$ is extended from $\overline D$ to $\overline D^*$ by setting $m(F)=0.$  An $m$-symmetric conservative diffusion process $X^*$ on $\overline D^*$ will be called a {\it symmetric conservative diffusion extension} of $X$ if its part process on $\overline D$ being killed upon hitting $F$ is equivalent in law with $X.$   The resolvent of $X$ is denoted by $\{G_\alpha^X,\; \alpha>0\}.$   

\begin{prop}\label{P:3.1} 
\ There exists a unique symmetric conservative diffusion extension 
$X^*$ of $X$ from $\overline D$ to $\overline D^*=\overline D\cup F$. 
The process $X^*$ is recurrent.  Let $(\EE^*, \FF^*)$ and $\FF^*_e$ be the Dirichlet form of $X^*$ on $L^2(\overline D^*,m)\;(=L^2(D;m))$ and its extended Dirichlet space, respectively.  Then
\begin{equation}\label{3.4}
\FF_e^*=H_e^1(D)\;\oplus\;\left\{\sum_{j=1}^N c_j\varphi_j\;:\; c_j\in \RR\right\}\;\subset\; {\rm BL}(D),
\end{equation}
\begin{equation}\label{3.5}
\EE^*(u,v)=\frac12\;\D(u,v),\qquad u,v \in \FF^*_e.
\end{equation}
\end{prop}

\pf  \ We apply a general existence theorem of a many-point extension 
formulated in \cite[Theorem 7.7.4]{CF2} to the $m$-symmetric diffusion $X$ on $\overline D$ and the $N$-points compactification $\overline D^*=\overline D\cup F$ of $\overline D.$  We verify conditions {\bf (M.1), (M.2), (M.3)} for $X$ required in this theorem.
$\psi_j(x):= \P_x(\zeta<\infty, X_{\zeta-}=\partial_j)$ is positive for every $x\in \overline D, 1\le j\le N,$ by \eqref{2.3} and \eqref{3.3}, and so {\bf(M.1)} is satisfied.  Since $m(\overline D)=\int_{\overline D}fdx<\infty,$ the $m$-integrability {\bf (M.2)} of the function 
$u_\alpha^{(j)}(x)=\E_x\left[e^{-\alpha\zeta};X_{\zeta-}=\partial_j\right],\ x\in \overline D,$ is trivially fulfilled,\ $1\le j\le N.$  For any $1\le j\le N$ and any compact set $V\subset \overline D,$ $\inf_{x\in V} G_\alpha^X\psi_j(x)$ is positive because  
$G_\alpha^X\psi_j=G^X_{0+} u_\alpha^{(j)}=G_{0+}^Z(u_\alpha^{(j)} f)$ 
is lower semi-continuous on account of \eqref{e:1.6} and $u_\alpha^{(j)}$ is positive on $\overline D.$   Accordingly, condition {\bf (M.3)} is also satisfied.

Therefore there exists an $m$-symmetric diffusion extension $X^*$ of $X$ from $\overline D$ to $\overline D^*$ admitting no killing on $F.$  We can then use a general characterization theorem \cite[Theorem 7.7.3]{CF2} to conclude that such an extension $X^*$ of $X$ is unique in law and its extended Dirichlet space $(\FF_e^*, \EE^*)$ is given by \eqref{3.4} and \eqref{3.5} as $\psi_j=\varphi_j,\; 1\le j\le N.$   In particular, \eqref{2.2} implies $1\in \FF_e^*,\ \EE^*(1,1)=0$, so that $X^*$ is recurrent and consequently conservative.  This also means the unique existence of an $m$-symmetric conservative diffusion extension $X^*$ of $X$ to $\overline D^*.$   
\qed

\begin{thm}\label{T:3.2}  ${\rm dim}(\HH^*(D))=N$ and
\begin{equation}\label{3.7}
\HH^*(D)=\left\{\sum_{j=1}^N c_j\; \varphi_j\;:\; c_j\in \RR\right\}.
\end{equation}
The $m$-symmetric conservative diffusion extension $X^*$ of the time changed RBM $X$ constructed in {\rm Proposition \ref{P:3.1}} is a reflecting extension of $X$ in the sense that the extended Dirichlet space $(\FF^*_e, \EE^*)$ of $X^*$ equals $({\rm BL}(D), \frac12\D)$ the reflected Dirichlet space of $X.$
\end{thm}

\pf\ By Proposition \ref{P:3.1}, 
$\{\varphi_j; 1\le j\le N\}\subset  \HH^*(D)\subset {\rm BL}(D).$
For $1\leq j, k\leq N$,  
let $c_k^{(j)} = c_k (\varphi_j )$. 
We claim that
\begin{equation}\label{3.8}
c_k^{(j)}=\delta_{jk},    \qquad 1\le k\le N.
\end{equation}
Let $\tau_n$ be the exit time of $Z$ from the set $\overline D\cap B_n(\o),\ n\ge 1.$  Then $\{\varphi_j(Z_{\tau_n})\}_{n\ge 1}$ is a bounded $\Q_x$-martingale
 and possesses an a.s. limit $\Phi$ with $\varphi_j(x)=\E^{Q_x}[\Phi].$
By \eqref{3.6},
\begin{equation}\label{3.9}
\Phi=\sum_{k=1}^N c_k^{(j)}\1_{\{Z_{\infty-}=\partial_k\}}.
\end{equation}
For $k\neq j,$ put $F_{k,n}=C_k\cap \{|x|=n\}.$  Then by \eqref{3.6} again
\begin{eqnarray*}
&&c^{(j)}_k\varphi_k(x)
=\lim_{n\to\infty}\E^{\Q_x}\left[\varphi_j(Z_{\tau_n})\1_{\{Z_{\infty-}=\partial_k\}}\right]
\le \limsup_{n\to\infty}\E^{\Q_x}\left[\varphi_j(Z_{\tau_n})\1_{\{Z_{\tau_n}\in C_k\}}\right]\\
&&=\limsup_{n\to\infty}
\E^{\Q_x}\left[\Q_x\left(Z_{\infty-}\circ\theta_{\tau_n}=\partial_j,\; Z_{\tau_n}\in C_k\big|\; \FF_{\tau_n}\right)\right]\\
&&\le\lim_{n\to\infty}
\Q_x\left(Z_{\infty-}=\partial_j,\; \sigma_{F_{k,n}}<\infty\right)=0,
\end{eqnarray*}
yielding $c_k^{(j)}=0,\ k\neq j.$  Taking $\Q_x$-expectation in \eqref{3.9} 
proves the  claim \eqref{3.8}.

Next for any $u\in {\rm BL}(D)$,   
let $ u_0= u-\sum_{j=1}^N c_j (u) \varphi_j$. 
Then $u_0\in {\rm BL}(D)$  with $c_j^{u_0}=0$
for every $1\leq j \leq N$. So by Proposition \ref{P:2.2}, $u_0\in H^1_e(D)$. 
This establishes \eqref{3.7}. The linear independence of $\{\varphi_j; 1\le j\le N\}$ follows 
from \eqref{3.8}, while \eqref{3.4} and \eqref{3.7} yield the last assertion of the theorem.  \qed

\begin{remark}\label{R:3.3} \rm This theorem for the special domain \eqref{e:1.5} was stated in \cite[Proposition 7.8.5]{CF2}.  We take this opportunity to mention that the proof of the latter given in the book \cite{CF2} contained a flaw (on the third line of page 386), that should be corrected in the above way. \qed 
\end{remark}

\section{Partitions $\Pi$ of $F$ and all possible symmetric diffusion extensions $Y$ of a time changed RBM $X$} \label{S:4}

We continue to consider the $N$-points compactification 
$\overline D^*=\overline D\cup F$  of $\overline D$ introduced at the end of Section 1.  A map $\Pi$ from the boundary set $F=\{\partial_1,\cdots, \partial_N\}$
onto a finite set $\wh F=\{\wh\partial_1,\cdots, \wh\partial_\ell\}$ with $\ell\le N$ is called a {\it partition} of $F.$  We let $\overline D^{\Pi,*}=\overline D\cup \wh F.$  We 
extend the map $\Pi$ from $F$ to $\overline D^*$ by setting $\Pi x=x,\ x\in \overline D,$ and introduce the quotient topology on $\overline D^{\Pi,*}$ by $\Pi$. 
In other words, we employ ${\cal U}_{\;\Pi}=\{U \subset \overline D^{\Pi,*}: \Pi^{-1}(U)\ \text{is an open subset of}\ \overline D^*\}$ as the family of open subsets of $\overline D^{\Pi,*}.$  
Then $\overline D^{\Pi,*}$ is a compact Hausdorff space and may be called an $\ell$-{\it points compactification of} $\overline D$ {\it obtained from} $\overline D^*$ {\it by identifying the points in the set} $\Pi^{-1}\wh \partial_i \subset F$ {\it as a single point} $\wh \partial_i$ {\it for each} $1\le i\le \ell.$ 

Given a partition $\Pi$ of $F$, the approaching probabilities $\wh\varphi_i$ of the RBM $Z=(Z_t, \Q_x)$ to $\wh \partial_i\in \wh F$ are defined by
\begin{equation}\label{4.1}
\wh\varphi_i(x)=\sum_{j\in \Pi^{-1}\wh\partial_i}\varphi_j(x),\quad x\in \overline D,\quad 1\le i\le \ell.
\end{equation}
As in the preceding section, we define the time changed process $X=(X_t,\zeta, \P_x)$ on $\overline D$ of $Z$ by means of a strictly positive bounded integrable function $f$ on $\overline D$.  The measure $m(dx)=f(x)dx$ is extended from $\overline D$ to $\overline D^{\Pi,*}$ by setting $m(\wh F)=0.$  Just as in 
Proposition \ref{P:3.1}, there exists then a unique $m$-symmetric conservative diffusion extension $X^{\Pi,*}$ of $X$ from $\overline D$ to $\overline D^{\Pi,*}$ and the Dirichlet form $(\EE^{\Pi,*}, \FF^{\Pi,*})$ of $X^{\Pi,*}$ on $L^2(\overline D^{\Pi,*};m)\;(=L^2(D;m))$ admits the extended Dirichlet space $(\FF_e^{\Pi,*}, \EE^{\Pi,*})$ expressed as
\begin{equation}\label{4.2}
\FF_e^{\Pi,*}=H_e^1(D)\;\oplus\;\left\{\sum_{i=1}^\ell c_i\wh\varphi_i\;:\; c_i\in \RR\right\}\;\subset\; {\rm BL}(D),
\end{equation}
\begin{equation}\label{4.3}
\EE^{\Pi,*}(u,v)=\frac12\;\D(u,v),\qquad u,v \in \FF^{\Pi,*}_e.
\end{equation}
$X^{\Pi,*}$ is recurrent.  $\EE^{\Pi,*}$ is a quasi-regular Dirichlet form on $L^2(\overline D^{\Pi,*};m).$\\

We now prove that the family $\{\overline X^{\Pi,*}: \Pi\ \text{is a partition of }\ F\}$ exhausts all possible $m$-symmetric conservative diffusion extensions of the time changed RBM $X$ on $\overline D$. 

Let $E$ be a Lusin space into which $\overline D$ is homeomorpically embedded as an open subset.
The measure $m(dx)=f(x)dx$ on $\overline D$ is extended to $E$ by setting $m(E\setminus \overline D)=0.$  Let $Y=(Y_t, \P_x^Y)$ be an $m$-symmetric conservative diffusion process on $E$ whose part process on $\overline D$ is identical in law with $X.$
 We denote by $(\EE^Y, \FF^Y)$ and $\FF^Y_e$ the Dirichlet form of $Y$ on $L^2(E;m)$ and its extended Dirichlet space.  
We call $Y$ an $m$-{\it symmetric conservative diffusion extension} of $X$. 
The following theorem extends \cite[Theorem 3.4]{CF1}.

\begin{thm}\label{T:4.1}
 There exists a partition $\Pi$ of $F$ such that, as Dirichlet forms  on $L^2(\overline D;m)$,
\begin{equation}\label{4.4}
(\EE^Y, \FF^Y)= (\EE^{\Pi,*}, \FF^{\Pi,*}).
\end{equation}
$Y$ under $\P_{g\cdot m}$ and $X^{\Pi.*}$ under $\P^{\Pi,*}_{g\cdot m}$ have the same finite dimensional distribution for any non-negative $g\in L^2(\overline D;m)$.  
Furthermore, a quasi-homeomorphic image of $Y$ is identical with $X^{\Pi,*}$ in the sense of 
{\rm Theorem \ref{T:6.2}} in Appendix. 
\end{thm}

\pf\ 
As has been noted in the preceding section, the extended Dirichlet space $(\FF_e^X, \EE^X)$ and the reflected Dirichlet space $((\FF^X)^{\rm ref},(\EE^X)^{\rm ref})$ of the Dirichlet form \eqref{3.2} are given by
\begin{equation}\label{4.5}
\FF_e^X=H_e^1(D),\qquad \EE^X=\frac12\D,
\end{equation}
\begin{equation}\label{4.6}
 (\FF^X)^{\rm ref}= {\rm BL}(D)=H_e^1(D)\oplus \HH^*(D),\qquad (\EE^X)^{\rm ref}=\frac12 \D,
\end{equation}
respectively.  

$\EE^Y$ is a quasi-regular Dirichlet form on $L^2(E;m)$ and $Y$ is properly associated with it by virtue of 
Z.-M. Ma and M. R\"ockner \cite{MR}.   By Chen-Ma-R\"ockner \cite{CMR}, 
$\EE^Y$ is therefore quasi homeomorphic with a regular Dirichlet form.   In particular, via a quasi homeomorphism $j$ in \cite[Theorems 3.1.13]{CF2}), we can  assume that $E$ is a locally compact separable metric space, $\EE^Y$ is a regular Dirichlet form on $L^2(E;m)$, $Y$ is an associated Hunt process on $E$, and $\wt F:=E\setminus \overline D$ is quasi-closed.  
Since $Y$ is a conservative extension of the non-conservative process $X$, $\wt F$ must be non $\EE^Y$-polar.  $Y$ can be also shown to be irreducible as in the proof of \cite[Lemma 7.2.7 (ii)]{CF2}.  
Thus we are in the same setting as in \S 7.1 of \cite{CF2} and Theorem 7.1.6 in it applies to $Y$ and $\wt F$.

Every function in $\FF^Y_e$ will be taken to be $\EE^Y$-quasi continuous.  
  As $Y$ is a diffusion with no killing inside, the jumping measure $J$ and the killing measure $k$ in the Beurling-Deny decomposition of $\EE^Y$ vanish so that we have by \cite[Theorem 7.1.6]{CF2}
\begin{equation}\label{4.7}
H_e^1(D)\subset \FF_e^Y\subset {\rm BL}(D),\quad \hh^Y:=\{\H u: u\in \FF_e^Y\}\subset \HH^*(D),
\end{equation}
\begin{equation}\label{4.8}
\EE^Y(u,u)=\frac12 \D(u,u)+\frac12 \mu_{\<\H u\>}^c(\wt F),\quad u\in \FF_e^Y,
\end{equation}  
where $\H u(x)=\E_x^Y[u(Y_{\sigma_{\wt F}})],\ x\in E.$

Let us prove that
\begin{equation}\label{4.9}
\mu_{\<u\>}^c(\wt F)=0\qquad u\in \hh^Y.
\end{equation} 
To this end, we consider a finite measure $\nu$ on $E$ defined by
\[\nu(B)=\int_{\overline D} \P_x^Y\left(Y_{\sigma_{\wt F}}\in B,\; \sigma_{\wt F}<\infty\right)m(dx),
\quad  B\in {\cal B}(E).
\]
$\nu$ vanishes off $\wt F$ and charges no $\EE^Y$-polar set.  In view of \cite[Lemma 5.2.9 (i)]{CF2}, $\wt F$ is a quasi support of $\nu$ in the following sense:
$\nu(E\setminus \wt F)=0$ and $\wt F\subset \wh F$ q.e. for any quasi closed set $\wh F$ with $\nu(E\setminus \wh F)=0.$

Now, for $u\in \hh^Y,$ \eqref{3.7} and \eqref{4.7} imply that $u=\sum_{j=1}^N c_j\varphi_j$ for some constants $c_j.$
Take $\wh F=\{\xi\in E: u(\xi)\in \{c_1,\cdots, c_N\}\}$.  
Since $u$ is quasi continuous, $\wh F$ is a quasi closed set.
As $u$ is continuous along the sample path of $Y$ (cf. \cite[Theorem 3.1.7]{CF2}), we have 
$\nu(E\setminus \wh F)=\P_m(u(Y_{\sigma_{\wt F}})\notin \{c_1,\cdots, c_N\})=0$  on account of 
Proposition \ref{P:2.2} and \eqref{3.8}.
Accordingly $\wt F\subset \wh F$ q.e., namely, 
 $u$ takes only finite values $\{c_1,\cdots, c_N\}$ q.e. on $\wt F$.
 By the {\it energy image density property} of $\mu^c_{\<u\>}$ 
due to N. Bouleau and F. Hirsch \cite{BH} 
(cf. \cite[Theorem 4.3.8]{CF2}), we thus get \eqref{4.9}.

Relation \eqref{4.7} and Proposition \ref{P:2.2}(ii) imply 
that every function $u\in \hh^Y (\subset {\rm BL}(D))$ admits a limit $u(\partial_j)$ at each boundary point $\partial_j\in F$ along the path of $Z$.
Define an equivalence relation $\sim$ on $F$ by $\partial_j\sim \partial_k$ if and only if $u(\partial_j)=u(\partial_k)$ for every $u\in \hh^Y.$  
Notice that, for every $1\le j\le N$, there exists $u\in \hh^Y$ with $u(\partial_j)\neq 0$.
Otherwise, for the resolvent $\{G_\alpha^Y: \alpha >0\}$ of $Y$, $G_\alpha^Y1\in \FF_e^Y (\subset {\rm BL}(D))$ approaches to zero at some $\partial_j$ along the path of $Z$, contradiction to the conservativeness of $Y$. 
Let $\Pi$ be the corresponding partition of $F$: $\Pi$ maps $F$ onto $\{\wh \partial_1,\cdots, \wh \partial_\ell\}$ the set of all equivalence classes with respect to $\sim.$  Then 
$\dis \hh^Y=\left\{\sum_{i=1}^\ell c_i\wh\varphi_i: \;c_i\in \RR \right\}$ for $\wh \varphi_i$ define by \eqref{4.1}.  
Hence \eqref{4.2}, \eqref{4.3}, \eqref{4.7}, \eqref{4.8} and \eqref{4.9} lead us to the desired identity \eqref{4.4}.

Since the both Dirichlet forms share a common semigroup on $L^2(\overline D;m)$, we get the first conclusion of the theorem.  Further the Dirichlet spaces
$$ (E, \; m, \; \EE^Y, \FF^Y),\qquad (\overline D^{\Pi,*},\; m,\; \EE^{\Pi,*}, \FF^{\Pi,*})
$$
are equivalent in the sense of Appendix (Section \ref{S:7}) by the identity map $\Phi$ from
 $\FF_b^Y$ onto $\FF_b^{\Pi,*}$ so that we get the second conclusion from Theorem \ref{T:6.2}.
\qed

\begin{remark}\label{R:4.2} \rm  
  (i)    For different choices of $f$, 
the family of all symmetric conservative extensions $Y$ of $X^f$ is invariant up to time changes because it shares a common family of extended Dirichlet spaces \eqref{4.2}-\eqref{4.3}. The same can be said for more general time changed RBM $X^\mu$,
which 
 will be formulated in Section \ref{S:6}. 

\smallskip

 (ii)   We can replace the conservativness assumption on $Y$ by a weaker one that $Y$ is a proper extension of $X$  
with no killing on $E\setminus \overline D.$
Then the above theorem remains valid if $X^{\Pi,*}$ is allowed to be replaced by its subprocess being killed 
upon hitting some (but not all) $\wh\partial_i$.
\qed 
\end{remark}

\begin{remark}\label{R:4.3}
\rm {\bf (Symmetric diffusion for a uniformly elliptic differential operator)}

 Given measurable functions $a_{ij}(x),\; 1\le i,j\le d,$ on $D$ such that
\begin{equation}\label{4.11}
a_{ij}(x)=a_{ji}(x),\quad \Lambda^{-1}|\xi|^2\le \sum_{1\le i,j\le d}a_{ij}(x)\xi_i\xi_j\le \Lambda|\xi|^2,\quad x\in D,\ \xi\in \RR^d,
\end{equation}
for some constant $\Lambda\ge 1,$ we consider a Dirichlet form 
\begin{equation}\label{4.12}
(\EE,\FF)=(\a, H^1(D))
\end{equation}
on $L^2(D)$ where
\[\a(u,v)=\int_D\sum_{i,j=1}^d a_{ij}(x)\frac{\partial u}{\partial x_i}(x)\frac{\partial v}{\partial x_j}(x)dx, \quad u,v\in H^1(D).\]
 If we replace the Dirichlet form \eqref{e:1.2} on $L^2(D)$ and the assoicated RBM $Z$ on $\overline D$, respectively, by the Dirichlet form \eqref{4.12} on $L^2(D)$ and the associated reflecting diffusion process on $\overline D$ constructed in \cite{FTo}, all results from Section \ref{S:2} to Section \ref{S:4} still hold without any  change as we shall see now.

By this replacement, the extended Dirichlet space and the reflected Dirichlet space are still $H_e^1(D)$ and BL$(D),$ respectively, although the inner product $\frac12\D$ is replaced by $\a.$  
The transience of \eqref{4.12} is equivalent to that of \eqref{e:1.2}.
The space $\HH^*(D)$ is now defined by \eqref{e:1.3} with $\a$ in place of $\frac12\D.$  
But, by noting that $\a(c,c)=0$ for any constant $c$ and by taking the characterization of a Liouville domain stated below Definition \ref{D:1.1} 
into account, 
we readily see that $D\in \DD$ is a Liouville domain relative to \eqref{4.12} if and only if so it is relative to \eqref{e:1.2}. \qed
\end{remark} 

\begin{remark}\label{R:4.4}  \rm {\bf (All possible symmetric conservative diffusion extensions of a one-dimensional minimal diffusion)}\quad
Consider a minimal diffusion $X$ on a one-dimensional open interval $I=(r_1,r_2)$ with no killing inside for which both boundaries $r_1, r_2$ are regular.  Let $E$ be a Lusin space into which $I$ is homeomorphically embedded as an open subset.  The speed measure $m$ of $X$ is extended to $E$ by setting $m(E\setminus I)=0.$  Let $Y$ be an $m$-symmetric conservative diffusion extension of $X$ from $I$ to $E.$  Then, by removing some 
$m$-polar open set for $Y$ from $\wt F=E\setminus I,$ a homeomorphic image of $Y$ is identical with 
either the two point extension of $X$ to $[r_1.r_2]$ or its one-point extension to the one-point compactification of $I.$     
This fact was implicitly indicated in \cite[\S 5]{F2} and \cite[\S 5]{F3} without proof.  This can be shown in a similar manner to the proof of Theorem \ref{T:4.1} 
by establishing the counterpart of the identity \eqref{4.9} and by noting that, for the one-point and two-point extensions of $X$, every non-empty subset of the state space has a positive $1$-capacity uniformly bounded away from zero due to  
 the bound \cite[(2.2.31)]{CF2} and so a quasi-homeomorphism is reduced to a homeomorphism.   

To put it another way, Theorem \ref{T:4.1} reveals that the time changed RBM $X$ on an unbounded domain with $N$-Liouville branches has a very similar structure to the one-dimensional diffusion only by changing two boundary points to $N$ boundary points. \qed
\end{remark} 

We note that the connected sum of non-parabolic manifolds being studied by Y. Kuz'menko and S. Molchanov \cite{KM}, A. Grigor'yan and L. Salloff-Coste \cite{GS} bears a strong similarity to the present paper in the setting although the main concern in these papers was the heat kernel estimates.  

\section{Characterization of $L^2$-generator of extension $Y$ by zero flux condition at infnity}
\label{S:5}

For a strictly positive bounded integrable function $f$ on $D$, we put $m(dx)=f(x)dx$ and denote by $(\cdot,\cdot)$ the inner product for $L^2(D;m).$
Let $Y$ be any $m-$symmetric conservative diffusion extension of the time changed process $X=X^f=(X_t, \zeta, \P_x)$ of the RBM $Z$ on $\overline D.$  Let 
$\Pi: F\mapsto \{\wh\partial_1,\cdots, \wh\partial_\ell\},\ \ell\le N,$ be the corresponding  partition of the boundary $F=\{\partial_1, \cdots, \partial_N\}$
appearing in Theorem \ref{T:4.1}.  The Dirichlet form $(\EE^Y,\FF^Y)$ of $Y$ on $L^2(D;m)$ is then described as
\[ 
\begin{cases}
\FF^Y= \left\{u = u_0+\sum_{i=1}^\ell c_i\wh\varphi_i:\; u_0\in H_e^1(D)\cap L^2(D;m),\ c_i\in \RR
\right\},\\
\EE^Y(u,v)=\frac12 \D(u,v),\qquad u, v\in \FF^Y,
\end{cases}
\]
where $\wh\varphi_i,\; 1\le i\le \ell,$ are defined by \eqref{4.1}.  

Let $\AA$ be the $L^2$-generator of $Y$, that is, $\AA$ is a self-adjoint operator on $L^2(D;m)$ such that $u\in \DD(\AA),\ \AA u=v\in L^2(D;m)$ if and only if $u\in \FF^Y$ with $\EE^Y(u,w)=-(v,w)$ for every $w\in \FF^Y.$  In view of Proposition \ref{P:2.2}, 
the condition (7.3.4) of \cite{CF2} is fulfilled by $Y.$  Therefore Theorem 7.7.3 (vii) of \cite{CF2} is well applicable in getting the following characterization of $\AA$:
\[u\in \DD(\AA)\quad \text{\rm if and only if}\quad u\in \DD(\LL)\ {\rm and}\ \NN(u)(\wh\partial_i)=0,\ 1\le i\le \ell.\] 
In this case, $\AA u=\LL u.$

Here $\LL$ is a linear operator defined as follows: $u\in \DD(\LL),\ \LL u=v\in L^2(D;m)$ if and only if $u\in {\rm BL}(D)\cap L^2(D;m)$ and $\frac12\D(u,w)=-(v,w)$ for every $w\in H_e^1(D)\cap L^2(D;m),$ or equivalently, for every $w\in C_c^1(\overline D).$\quad
$\NN(u)(\wh\partial_i)$ is the {\it flux of} $u$ {\it at} $\wh\partial_i$ defined by
\[\NN(u)(\wh\partial_i)=\frac12\D(u,\wh\varphi_i)+(\LL u,\wh\varphi_i),\quad 1\le i\le \ell.\]

It can be readily verified that $u\in \DD(\LL)$ if and only if $u\in {\rm BL}(D)\cap L^2(D;m),$ $\Delta u$ in the Schwartz distribution sense is in $L^2(D)$ and
\begin{equation}\label{5.1}
\D(u,w)\;+\; \int_D\Delta u(x)\cdot w(x) dx\; =\; 0\quad\text{\rm for every}\ w\in C_c^1(\overline D).
\end{equation} 
In this case, 
$\LL u(x)=\frac{1}{2f(x)}\;\Delta u(x),\ x\in D.$
The equation \eqref{5.1} can be interpreted as the requirement that the {\it generalized normal derivative} 
of $u$ vanishes on $\partial D.$
Thus we have

\begin{thm}\label{T:5.1}  $u\in \DD(\AA)$ if and only if $u\in {\rm BL}(D)\cap L^2(D;m)$, $\Delta u$ in the Schwartz distribution sense belongs to $L^2(D),$ the equation {\rm \eqref{5.1}} is satisfied and
\begin{equation}\label{5.2}
\left(\NN(u)(\wh\partial_i)=\right) \frac12\D(u,\wh\varphi_i)+\frac12\int_D\Delta u(x)\wh\varphi_i(x)dx=0,\quad 1\le i\le \ell.
\end{equation}
In this case,
\begin{equation}\label{5.3}
\AA u(x)=\frac{1}{2f(x)}\;\Delta u(x),\quad\text{\rm a.e. on}\ D.
\end{equation}
\end{thm}

Suppose $u\in \DD( \AA)$ is smooth on $\overline D.$  Then $\frac{\partial u}{\partial \n}=0$ on $\partial D$ due to the condition \eqref{5.1} so that the zero flux condition \eqref{5.2}  
at $\wh \partial_j$ can be expressed as
\begin{equation}\label{5.4}
\lim_{r\uparrow \infty}\int_{D\cap \partial B_r(\o)} u_r(x) \wh \varphi_i(x)d\sigma_r(dx)=0,\quad 1\le i\le \ell,
\end{equation} 
where $\sigma_r$ is the surface measure on $\partial B_r(\o).$

The last part of Section 7.6 $(4^\circ)$ of \cite{CF2} has treated a very special case of the above where $D=\RR^d,\: d\ge 3,$ and $Y$ is the one-point reflection at the infinity of $\RR^d$ of a time changed Brownian motion on $\RR^d.$

In \cite{F3}, the $L^2$-generator of any symmetric diffusion extension $Y$ of a one-dimensional minimal diffusion $X$ is identified.  In this case, the Dirichlet form of $Y$ admits its reproducing kernel which enables us to identify also the $C_b$-generator of  $Y$, 
 recovering the general boundary condition due to 
W. Feller and K. It\^o-H. P. McKean.

\section{Extensions of more general time changed RBMs} \label{S:6}

All the results in Sections \ref{S:3}-\ref{S:5} except for \eqref{5.3} hold  for more general time changed RBMs than $X^f$.
Let $Z=(Z_t, \Q_x)$, $f$, $X=X^f=(X_t,\zeta,\P_x)$, $X^*=(X_t^*, \P_x^*)$ be as in Section \ref{S:3}. 

We consider a finite smooth measure $\mu$ on $\overline D$ with full quasi-support $\overline D$ relative to the Dirichlet form $(\EE,\FF)$ of \eqref{e:1.2}.  Let $A^\mu$ be the PCAF of $Z$ with Revuz measure $\mu$ and $X^\mu=(X_t^\mu, \zeta^\mu, \P_x^\mu)$ be the time changed process of $Z$ by $A^\mu$.
The Markov process   $X^\mu$ is $\mu$-symmetric and its Dirichlet form 
$(\EE^{X^\mu}, \FF^{X^\mu})$ 
on $L^2(\overline D;\mu)$ is given by
\begin{equation}\label{eq:6.1}
\EE^{X^\mu}=\frac12\D,\quad \FF^{X^\mu}=H_e^1(D)\cap L^2(\overline D;\mu). 
\end{equation}

\begin{prop}\label{P:6.1}  It holds that
\begin{equation}\label{eq:6.2}
\Q_x(A_\infty^\mu<\infty)=1\quad
\hbox{\rm for q.e. } 
x\in \overline D,
\end{equation}
 \begin{equation}\label{eq:6.3}
\P_x^\mu(\zeta^\mu<\infty, \ X_{\zeta^\mu-}^\mu =\partial_i)=\varphi_i(x)>0,\quad\hbox{\rm for q.e. } 
x\in \overline D \hbox{ \rm and }  1\le i \le N.
\end{equation} 
\end{prop}

\pf  Fix a strictly positive bounded integrable function $h_0$.
By the transience of $Z$ and   
\cite[Theorem A.2.13 (v)]{CF2}, $G_{0+}^Z h_0(x)<\infty$ 
for q.e. $x\in \overline D$.
For integer $k\geq 1$,  
let 
$$  
\Lambda_k:= \left\{ x\in \overline D: G_{0+}^Z h_0(x)\le 2^k \right\} \quad 
\hbox{and} \quad  
 h(x)= \sum_{k=1}^\infty 2^{-2k} {\bf 1} _{\Lambda_k }(x) h_0(x).
 $$
 Then $h$ is a strictly positive bounded integrable function on $\overline D$ with  
   $G_{0+}^Z h(x)\le \1$ q.e. on $\overline D$.
  From \cite[(4.1.3)]{CF2}, we have
\begin{equation}\label{eq:6.4}
\int_{\overline D} \E^{\Q_x}\left[A_\infty^\mu\right]h(x)dx=\< G_{0+}^Z h,\mu\>\leq \mu(\overline D)<\infty.
\end{equation}
  It follows that $\E^{\Q_x}\left[A_\infty^\mu\right]<\infty$   
a.e $x\in \overline D$ and hence q.e. $x\in \overline D$ by
\cite[Theorem A.2.13 (v)]{CF2}, yielding
\eqref{eq:6.2}.   \eqref{eq:6.3} follows from \eqref{eq:6.2} and Proposition \ref{P:2.1}. 
\qed

Since $m(dx)=f(x)dx$ has its quasi-support $\overline D$ relative to $(\EE,\FF)$, the Dirichlet form $(\EE^X,\FF^X)$ of 
\eqref{3.2}
shares the common quasi-notion with $(\EE,\FF)$ (\cite[Theorem 5.2.11]{CF2}).  Hence the quasi-support of $\mu$ relative to $(\EE^X,\FF^X)$ is still $\overline D$.

The Dirichlet form $(\EE^*, \FF^*)$ on $L^2(\overline D^*,m)$ 
of $X^*$ is quasi-regular.  According to the quasi-homeomorphism method 
already used in Section 4, we may assume it to be regular.  The measure $\mu$ on $\overline D$ is extended to $\overline D^*$ by setting $\mu(F)=0$.  We claim that the quasi-support of $\mu$ relative to this Dirichlet form equals $\overline D^*$ by using a criteria \cite[Theorem 3.3.5]{CF2}.

Assume that $u\in \FF^*$ is $\EE^*$-quasi-continuous and that $u=0$  $\mu$-a.e.
Then $u\big|_{\overline D}$ is $\EE^X$-quasi-continuous (\cite[Theorem 3.3.8]{CF2}) so that $u=0$ q.e. on $\overline D$.  According to the same reference, there exists a Borel $m$-polar set $C\subset \overline D$ relative to $X^*$ such that $u(x)=0$ for every $x\in \overline D\setminus C.$ 
Since $u$ is continuous along the path of $X^*$ (\cite[Theorem 3.1.7]{CF2}), we have for each $1\le i\le N$
\[ \P_m^*\left( u(\partial_i)=\lim_{t\uparrow \sigma_F} u(X_t^*),\; \sigma_C=\infty,\; \sigma_F<\infty,\;  X_{\sigma_F}^*=\partial_i\right)=\P_m(\zeta<\infty, X_{\zeta-}=\partial_i)>0,\]
and so $u$ vanishes on $F$ and hence q.e. on $\overline D^*$, as was to be proved.

\begin{thm}\label{thm:6.2}\ There exists  a unique 
 $\mu$-symmetric conservative diffusion $\wt X^{*,\mu}$ on $\overline D^*$ which is a q.e. 
extension of $X^\mu$ in the sense that the part of the former on $\overline D$ coincides in law with the latter for 
q.e.  starting points $x\in \overline D.$   The extended Dirichlet space of $\wt X^{*,\mu}$ equals 
$({\rm BL}(D), \frac12\D)$ the reflected Dirichlet space of $X^\mu.$  
\end{thm}

\pf\ Let $B_t^0$ and $B_t$ be the PCAFs of $X$ and $X^*$, respectively, with Revuz measure $\mu.$  According to \cite[Proposition 4.1.10]{CF2}
\begin{equation}\label{eq:6.5}
B_t^0=B_{t\wedge\sigma_F}.
\end{equation}
Let $\wt X^\mu$ and $ \wt X^{*,\mu}$ be the time changed processes of $X$ and $X^*$ by means of $B_t^0$ and $B_t$, respectively.   
The Markov process $\wt X^\mu$ is then the part of $\wt X^{*,\mu}$ on $\overline D$ by \eqref{eq:6.5}.  Since $X^*$ is recurrent, so is $\wt X^{*,\mu}$ in view of \cite[Theorem 5.2.5]{CF2}.  Therefore $\wt X^{*,\mu}$ is a $\mu$-symmetric conservative diffusion extension of $\wt X^\mu.$

On the other hand, 
the Dirichlet form of $\wt X^\mu$ on $L^2(\overline D;\mu)$ 
is identical with \eqref{eq:6.1}    
the Dirichlet form of $X^\mu$ on $L^2(\overline D;\mu)$, and consequently $\wt X^{*,\mu}$ is a q.e. extension of $X^\mu.$ 
The last statement follows from the invariance of extended and reflected Dirichlet spaces under time changes by fully supported PCAFs. 

The uniqueness of such a $\mu$-symmetric conservative Markovian extension of $X^\mu$  to $\overline D^*$ follows from   \cite[Theorem 7.7.3]{CF2} .
\qed 

Similarly, all results in Section 4 and 5 with $\mu$ in place of $dm=fdx$ remain valid except for \eqref{5.3}.

\begin{remark}\label{rem:6.3}\rm\ One can give an alternative proof of  Theorem \ref{thm:6.2} 
without invoking the time change of $X^*$ but still using the quasi-regularity of $(\EE^*,\FF^*)$.  Indeed, the following proposition combined with \eqref{eq:6.3} and \cite[Theorem 7.7.3]{CF2} readily yields Theorem \ref{thm:6.2}.
\end{remark}

Each function in $\FF^*_e$ is taken to be $\EE^*$-quasi continuous.  Define 
\begin{equation}\label{eq:6.6}
\wh \FF =\FF^*_e\cap L^2(\overline D; \mu)  
\quad \hbox{\rm and} \quad 
\wh\EE (u, v)=\EE^*(u, v)= \frac12 \D (u, v) \ \hbox{\rm for } u, v \in \wh \FF.
\end{equation}

\begin{prop}\label{prop:6.4}\begin{description}
\item {\rm (i)}\ $( \wh \EE, \wh \FF)$ is a quasi-regular Dirichlet form on $L^2(\overline D^*; \mu)$.   

 \item{\rm (ii)} 
Its  associated strong Markov process 
$\wh X$ on $\overline D^*$ is a 
 $\mu$-symmetric conservative diffusion 
which is a q.e. extension of $X^\mu$.

  \item{\rm (iii)} Each $\partial_j$ is non-$\wh \EE$-polar. 
\end{description}
\end{prop}

\pf\ (i)\ As $\overline D$ is a quasi-support of $\mu$, $u=0\; \mu-$a.e. for $u\in \wh\FF$ implies $u=0$ a.e. on $\overline D$ and $\D(u,u)=0$.
This together with the transience of $(\FF^*_e, \EE^*)$ implies that 
 $( \wh \EE, \wh \FF)$ is a well defined Dirichlet form on
$L^2(\overline D^*; \mu)$.

Since $(\EE^*,   \FF^* )$ is a quasi-regular Dirichlet form on $L^2(\overline D^*; m)$,
by \cite[Remark 1.3.9]{CF2},  there is an increasing sequence of 
compact subsets $\{F_k\}$ of $\overline D^*$ so that 
\begin{description}
\item{(a)} there is an increasing sequence of 
compact subsets $\{F_k\}$ of $\overline D^*$ so that $\cup_{k\geq 1}  \FF^*_{F_k}$
is $\EE^*_1$-dense in $\FF^*$.
\item{(b)} there is an $\EE^*_1$-dense of  countable set $\Lambda_0:=\{f_j; j\geq 1\}$ of bounded functions of $\FF^*$
so that $ \{f_j; j\geq 1\}\subset C(\{F_k\})$ and they separate points of $\cup_{k\geq 1} F_k$. 
\end{description}
By the contraction of the Dirichlet form, we may and do assume without loss of generality that 
for every integer $n\geq 1$ and $f\in \Lambda_0$, $((-n)\vee f)\wedge n \in \Lambda_0$. 
We claim that $\cup_{k\geq 1}  \FF^*_{F_k, b}\subset \cup_{k\geq 1}  \wh \FF_{F_k, b}$
is $\wh \EE_1$-dense in $\wh \FF_b$.  
Let $u\in \wh \FF_b$. Since $\wh \FF_b =\FF^*_b$, there are $u_k\in \FF^*_{F_k}$
so that $u_k\to u$ in $\EE^*_1$-norm.
Using truncation if needed, we may and do assume  $\|u_k\|_\infty \leq \|u\|_\infty +1$. 
 Taking a subsequence if needed, we may also assume that $u_k$ converges to $u$ $\EE^*$-q.e. on  $\overline D^*$.
 Since $\mu$ is a finite smooth measure, we conclude that $u_k$ is $\wh \EE_1$-convergent to $u$.
 This proves the claim. As $\wh \FF_b$ is $\wh \EE_1$ dense in $\wh \FF$,
 it follows that  $\{F_k\}$ is an $\wh \EE$-nest on $\overline D^*$.
 
 A similar argument shows that $\Lambda_0\subset  \wh \FF_b=\FF^*_b$ is $\wh \EE_1$-dense in
 $\wh \FF_b$ and hence in $\wh \FF$.  This proves the assertion (i).
 
\noindent
(ii)\ Since $1 \in \wh \FF$ and $\D(1,1)=0$, the associated  $\mu$-symmetric diffusion $\wh X$ on 
$\overline D^*$ is recurrent and conservative. 
For $R>r$, take $\psi\in C_c^\infty(\overline D)$ with $\psi=1$ on $B_{R+1}({\bf 0})$.   Then, for any bounded $u\in \wh \FF$, $\psi u\in H_e^1(D)$ and so 
\[\{v\in \wh \FF: v=0\ \hbox{\rm q.e. on } \overline D^*\setminus B_R({\bf 0})\}
=\{v\in H_e^1(D)\cap L^2(\overline D;\mu): v=0\ \hbox{\rm q.e. on } \overline D\setminus B_R({\bf 0})\},\]
 namely, the part of $\wh\EE$ on $\overline D \cap B_R({\bf 0})$ coincides with  the part of $\EE^{X^\mu}$ on $\overline D \cap B_R({\bf 0}).$
By letting $R\to\infty$, we see that 
the part of $\wh\EE$ on $\overline D$ coincides with  $\EE^{X^\mu},$ proving (ii).    

 \noindent (iii) The non-$\wh \EE$-polarity 
of $\partial_j$ follows from (ii) and \eqref{eq:6.3}. 
\qed

\section{Appendix: equivalence and quasi-homeomorphism}\label{S:7}

In dealing with boundary problems for symmetric Markov processes,
 it is convenient to introduce an equivalence of Dirichlet spaces following \cite[A.4]{FOT} as will be stated below.
 
We say that a quadruplet $(E,m,\EE,\FF)$ is a {\it Dirichlet space} if $E$ is a Hausdorff topological space with a countable base, $m$ is a $\sigma$-finite positive Borel measure on $E$ and $\EE$ with domain $\FF$ is a Dirichlet form on $L^2(E;m)$. The inner product in $L^2(E;m)$ is denoted by $(\cdot,\cdot)_E$.
For a given Dirichlet space $(E, m, \EE, \FF)$, the notions of an $\EE$-{\it nest}, an $\EE$-{\it polar set}, an $\EE$-quasi-{\it continuous numercal function} and \lq $\EE$-{\it quasi-everywhere}' (\lq $\EE$-q.e.' in abbreviation) are defined as in \cite[Definition 1.2.12]{CF2}.  The {\it quasi-regularity} of the Dirichlet space is defined just as in \cite[Definition 1.3.8]{CF2}.
We note that the space $\FF_b= \FF \cap L^\infty(E;m)$ is an algebra.

\begin{remark}\label{R:6.1} \rm  In Section 1.2 and the first half of Section 1.3 of \cite{CF2}, it is assumed that 
\begin{equation}\label{6.1}
{\rm supp}[m]=E.  
\end{equation}
We need not assume it.  Generally, if we let $E'={\rm supp}[m]$, then $E\setminus E'$ is $\EE$-polar 
according to the definition of the $\EE$-polarity.
  If $(E,m, \EE,\FF)$ is quasi-regular, so is $(E', m\big|_{E'},\EE,\FF)$ accordingly.  Therefore we may assume \eqref{6.1} if we like by replacing $E$ with $E'$. \qed 
\end{remark}

Given two Dirichlet spaces
\begin{equation}\label{6.2}
 (E,m,\EE,\FF),\qquad (\tilde{E},\tilde{m},\tilde{\EE},\tilde{\FF}), 
\end{equation}
we call them {\it equivalent} if there is an algebraic isomorphism $\Phi$ from $\FF_b$ onto $\tilde{\FF}_b$ preserving three kinds of metrics: for $u \in \FF_b$
\[\Vert u\Vert_\infty = \Vert\Phi u\Vert_\infty,\ (u,u)_E = (\Phi u,\Phi u)_{\tilde{E}},\ \EE(u,u) = \tilde{\EE}(\Phi u, \Phi u).
\]
One of the two equivalent Dirichlet spaces is called a {\it representation} of the other.

The underlying spaces $E, \ \tilde{E}$ of two Dirichlet spaces \eqref{6.2} are said to be {\it quasi-homeomorphic} if there exist $\EE$-nest $\{F_n\}$, $\tilde{\EE}$-nest $\{\tilde{F}_n\}$ and a one to one mapping $q$ from $E_0 = \cup_{n=1}^\infty F_n$ onto $\tilde{E}_0 = \cup_{n=1}^\infty \tilde{F}_n$ such that the restricition of $q$ to each $F_n$ is a homeomorphism onto $\tilde F_n$.  $\{F_n\},\ \{\tilde F_n\}$ are called the {\it nests attached to the quasi-homemorphism} $q$. 
Any quasi-homeomorphism is quasi-notion-preseving.

We say that the equivalnce $\Phi$ of two Dirichlet spaces \eqref{6.2} {\it is induced by a quasi-homeomorphism} $q$ of the underlying spaces if
\[ \Phi u(\tilde x)=u(q^{-1}(\tilde x)),\quad u\in \FF_b,\quad \tilde m{\rm -a.e.}\ \tilde x.\]
Then $\tilde m$ is the image measure of $m$ and $(\tilde \EE, \tilde \FF)$ is the image Dirichlet form of $(\EE,\FF).$  

\begin{thm}\label{T:6.2} Assume that two Dirichlet spaces {\rm\eqref{6.2}} are quasi-regular and that they are equivalent.  Let $X=(X_t, \PP_x)$ {\rm(}resp. $\tilde X=(\tilde X_t, \tilde \PP_x)${\rm)} be 
an $m$-symmetric right process on $E$ 
{\rm(}resp. an $\tilde m$-symmetric right process on $\tilde E${\rm)} properly associated with 
$(\EE,\FF)$ on $L^2(E;m)$ 
{\rm(}resp. $(\tilde\EE,\tilde\FF)$ on $L^2(\tilde E;\tilde m)${\rm)}.
Then the equivalence is induced by a quasi-homeomorphism $q$ with attached nests $\{F_n\},$ $\{\tilde F_n\}$ such that $\tilde X$ is the image of $X$ by $q$ in the following sense: there exist an $m$-inessential Borel subset 
$N$ of $E$ containing $\cap_{n=1}^\infty F_n^c$ 
and an $\tilde m$-inessential Borel subset 
$\tilde N$ of $\tilde E$ containing $\cap_{n=1}^\infty \tilde F_n^c$ so that $q$ is one to one from $E\setminus N$ onto $\tilde E\setminus \tilde N$ and
\begin{equation}
\tilde X_t=q(X_t),\qquad \tilde \PP_{\tilde x}=\PP_{q^{-1}\tilde x},\qquad \tilde x\in \tilde E\setminus \tilde N.
\end{equation} 
\end{thm}

\pf\ Since both Dirichlet spaces in \eqref{6.2} are assumed to be quasi-regular, they are equivalent to some regular Dirichlet spaces and the equivalences are induced by some quasi-homeomorphisms $q_1, q_2$ in view of \cite[Theorem 1.4.3]{CF2}.  Since two Dirichlet spaces in \eqref{6.2} are also assumed to be equivalent, so are the corresponding two regular Dirichlet spaces, the equivalence being induced by a quasi-homeomorphism $q_3$ on account of \cite[Theorem A.4.2]{FOT} combined with \cite[Theorem 1.2.14]{CF2}.    Hence the equivalence of the quasi-regular Dirichlet spaces in \eqref{6.2} is induced by the quasi-homeomorphism $q=q_1\circ q_3\circ q_2^{-1}$ between $E$ and $\tilde E$.  Let $\{F_n\}, \{\tilde F_n\}$ be the nests attached to $q$.

According to \cite[Theorem 3.1.13]{CF2}, we may assume without loss of generality that both $X$ and $\tilde X$ are Borel right processes.  Further the $\EE$-polarity is equivalent to the $m$-polar for $X.$  By virtue of \cite[Theorem A.2.15]{CF2}, we can therefore find an $m$-inessential Borel set $N_1\subset E$ containing $\cap_{n=1}^\infty F_n^c.$  Consider the set $\tilde N_1\subset \tilde E$ defined by $q(E\setminus N_1)=\tilde E\setminus \tilde N_1.$  $\tilde N_1$ is an $\tilde \EE$-polar Borel set and $q$ is one to one from $E\setminus N_1$ onto $\tilde E\setminus \tilde N_1.$

Define the process $\wh X=(\wh X_t, \wh\PP_{\tilde x})_{\tilde x\in \tilde E\setminus \tilde N_1}$ by
\[
\wh X_t=q(X_t),\qquad \wh \PP_{\tilde x}=\PP_{q^{-1}\tilde x},\qquad \tilde x\in \tilde E\setminus \tilde N_1.
\] 
On account of \cite[Lemma 3.1]{FFY}, we can then see that $\wh X$ is an $\tilde m$-symmetric Markov process on $\tilde E\setminus \tilde N_1$ 
properly associated with the Dirichlet form $(\tilde \EE,\tilde \FF)$ on $L^2(\tilde E;\tilde m).$
  Since the $\tilde m$-symmetric Borel right process $\tilde X$ is also
properly associated with the Dirichlet form $(\tilde \EE,\tilde \FF)$ on $L^2(\tilde E;\tilde m),$ the same method as in the proof of \cite[Theorem 3.1.12]{CF2} combined with \cite[Theorem A.2.15]{CF2} leads us to finding an $\tilde m$-inessential Borel set $\tilde N$ containing $\tilde N_1$ for $\tilde X$ such that 
the Markov processes $\tilde X\big|_{\tilde E\setminus \tilde N}$  
 and $\wh X\big|_{\tilde E\setminus \tilde N}$ are identical in law.
It now suffices to define the set $N$ by $E\setminus N=q^{-1}(\tilde E\setminus \tilde N).$  \qed  

\begin{remark}\label{R:6.3} \rm 
Owing to the works of S. Albeverio, Z.-M. Ma, M. R\"ockner and P. J. Fitzsimmons, the quasi-regularity of a Dirichlet form has been known to be not only a sufficient condition but also a necessary one for the existence of a properly associated right process.   
It is further shown in \cite{CMR} that a Dirichlet form is quasi-regular if and only if it is quasi-homeomorphic to a reglar Dirichlet form on a locally compact separable metric space. 
These facts are formulated by Theorem 1.5.3 and Theorem 1.4.3, respectively, of \cite{CF2}
under the assumption \eqref{6.1} which is not needed actually.
But we may assume it without loss of generality as will be seen below.

Indeed, let $E$ be a Lusin space, $m$ be a $\sigma$-finite measure on $E$ and $X$ be an $m$-symmetric Borel right process on $E.$  
Then, for $E_0={\rm supp}[m]$, $E\setminus E_0$
 is an $m$-negligible open set so that it is $m$-polar for $X$ by \cite[Theorem A.2.13 (iii)]{CF2}.  Hence, by \cite[Theorem A.2.15]{CF2}, there exists a Borel set $E_1\subset E_0$ such that $E\setminus E_1$ is $m$-inessential for $X$.
$E_1$ is the support of $m\big|_{E_1}$ because, for any $x\in E_1$ an any neighborhood $O(x)$ of $x$, $m(O(x)\cap E_1)=m(O(x))-m(O(x)\cap (E\setminus E_1))>0$.
Hence it suffices to replace $E$ by $E_1$. 

In Theorem \ref{T:4.1}, the extension process $Y$ is assumed to live on a Lusin space $E$ into which $\overline D$ is homeomorphically embedded as an open subset. In this particular case, the above set $E_1$ can be choosen to contain $\overline D$ on account of the proof of \cite[Theorem A.2.15]{CF2}.  Therefore, in Theorem \ref{T:4.1} (resp. Remark \ref{R:4.4}), we can assume more strongly that $\overline D$ (resp. $I$)  is homeomorphically embedded into the state space $E$ of $Y$ as a dense open subset. \qed 
\end{remark}

{\small

\vskip 0.3truein

\noindent {\bf Zhen-Qing Chen}

\smallskip \noindent
Department of Mathematics, University of Washington, Seattle,
WA 98195, USA

\noindent
E-mail: \texttt{zqchen@uw.edu}

\medskip

\noindent {\bf Masatoshi Fukushima}:

\smallskip \noindent
 Branch of Mathematical Science,
Osaka University, Toyonaka, Osaka 560-0043, Japan.

\noindent Email: \texttt{fuku2@mx5.canvas.ne.jp}

\end{document}